\documentclass[reqno]{amsart}
\usepackage{latexsym,amsfonts,amssymb,amsmath,euscript}
\usepackage[latin1]{inputenc} 
\usepackage{float}
\usepackage{graphicx} % inserir imagens
\usepackage{color,fancybox
}
\newcommand{\R} {\mathbb{R}}
\newcommand{\N} {\mathbb{N}}

\newcommand{\eps}{\epsilon}

\begin{document}
\title{Thin domains with doubly oscillatory boundary $^\dag$}
%\author[Arrieta and Pereira]
\author[J. M. Arrieta]{Jos\'{e} M. Arrieta}
%\thanks{$^*$ Corresponding author:  Jos\'e M. Arrieta,  Departamento de Matem\'atica Aplicada, Facultad de Matem\'aticas, 
%Universidad Complutense de Madrid, 28040  Madrid, Spain. e-mail: arrieta@mat.ucm.es}

\address[Jos\'e M. Arrieta]{Departamento de Matem\'atica Aplicada,
Facultad de Ma\-te\-m\'a\-ti\-cas, Universidad Complutense de
Madrid, 28040 Madrid, Spain.} \email{arrieta@mat.ucm.es}

\author[M. Villanueva-Pesqueira]{Manuel Villanueva-Pesqueira.}
\thanks{$^\dag$ Partially supported by grant MTM2012-31298, MINECO, Spain and Grupo de Investigaci\'on CADEDIF, UCM}
\address[M. Villanueva-Pesqueira]{Departamento de Matem\'atica Aplicada,
Facultad de Ma\-te\-m\'a\-ti\-cas, Universidad Complutense de
Madrid, 28040 Madrid, Spain.} \email{manuelvillanueva@mat.ucm.es}

\date{}

\begin{abstract}

We consider a 2-dimensional thin domain with order of thickness $\eps$ which presents oscillations of amplitude also $\epsilon$ on both boundaries , top and bottom, but the period of the oscillations are of different order at the top and at the bottom. We study the behavior of the Laplace operator with Neumann boundary condition and obtain its asymptotic homogenized limit as $\eps\to0$. We are interested in understanding how this different oscillatory behavior at the boundary, influences the limit problem.  
\end{abstract}

\maketitle
\numberwithin{equation}{section}
\newtheorem{theorem}{Theorem}[section]
\newtheorem{lemma}[theorem]{Lemma}
\newtheorem{corollary}[theorem]{Corollary}
\newtheorem{proposition}[theorem]{Proposition}
\newtheorem{remark}[theorem]{Remark}
\allowdisplaybreaks

\section{Introduction} \label{introduction}
\begin{figure}[H]
  \centering
    \includegraphics [width=0.3\textwidth]{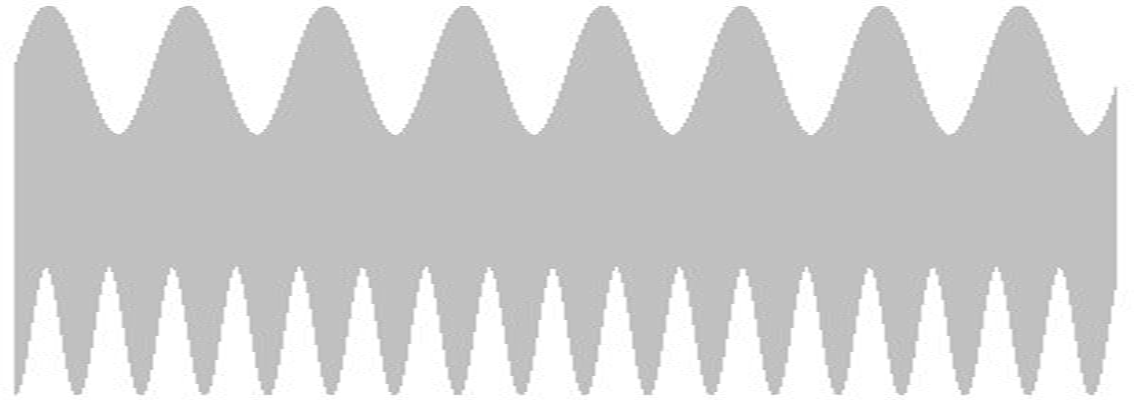}
    \caption{Thin domain $R^\eps$}
    \label{thin}
\end{figure}

In this paper, we analyze the behavior of the solutions of the Laplace equation with homogeneous Neumann boundary conditions
\begin{equation} \label{OPI0}
\left\{
\begin{gathered}
- \Delta w^\epsilon + w^\epsilon = f^\epsilon
\quad \textrm{ in } R^\epsilon \\
\frac{\partial w^\epsilon}{\partial N^\epsilon} = 0
\quad \textrm{ on } \partial R^\epsilon
\end{gathered}
\right. 
\end{equation}
with $f^\epsilon \in L^2(R^\epsilon)$ and $N^\epsilon$ is the unit outward normal to $\partial R^\epsilon$. 
The domain $R^\epsilon$ is a  two dimensional thin domain which presents a highly oscillatory behavior at the boundary and it is given as the region between two oscillatory functions, that is, 
\begin{equation}\label{thin-intro}
R^\epsilon = \Big\{ (x_1,x_2) \in \R^2 \; | \;  x_1 \in (0,1),  \;
 -\eps \,h(x_1/\eps^\alpha) < x_2 < \epsilon \, g(x_1/\eps) \Big\},  \quad \text{with} \,\, \alpha>1.
 \end{equation}
where  $g, h : \R \to \R$ are $C^1$ periodic functions with period $L_1$ and $L_2$ respectively (see Figure \ref{thin}). Moreover, there exist constants $h_0\geq 0$ and  $h_1,g_0,g_1>0$ such that 
$
0\leq h_0\leq h(\cdot)\leq h_1,
$
and
$
0<g_0 \leq g(\cdot)\leq g_1.
$  

Observe that both the amplitude and period of the oscillations at the upper boundary, given by $\eps g(x/\eps)$  are of the same order as the thickness of the domain. But, for the lower boundary, which is given by $\eps h(x/\eps^\alpha)$, the amplitude is of the same order $\eps$, while the period is of the order of $\eps^\alpha$, which means that we have much more oscillations at the bottom than at the top boundary. order of the lower oscillations is large than the order of the amplitude and height of the thin domain $R^\eps$ with respect to the small parameter 
$\eps$.

The existence and uniqueness of solutions for problem (\ref{OPI0}) for each $\eps>0$, is guaranteed by Lax-Milgram Theorem. We will analyze the asymptotic behavior of  the solutions as $\eps \to 0$.

Since the domain is thin,  $R^\epsilon \subset (0,1) \times (-\eps h(\cdot), \eps g(\cdot))$, approaching the interval $(0,1)$, it is reasonable to expect that
the family of solutions will converge to a function of just one variable and that this function will satisfy
certain elliptic equation  in one dimension with some boundary conditions. As a matter of  fact, if the function $h_\eps(\cdot)$ is independent of $\eps$, say $h_\eps(\cdot)\equiv 0$,  the limit equation is given by
\begin{equation} \label{PC}
\left\{
\begin{gathered}
-q_0w_{xx} + w =  f(x), \quad x \in (0,1)\\
w'(0)=w'(1)=0
\end{gathered}
\right.
\end{equation} 
where
$
q_0=\frac{1}{|Y^*|} \int_{Y^*} \Big\{ 1 - \frac{\partial X}{\partial y_1}(y_1,y_2) \Big\} dy_1 dy_2,
$ 
and $X$ is a convenient auxiliary harmonic function
defined in the representative basic cell
$
Y^*= \{ (y_1,y_2) \in \R^2 \; | \; 0< y_1 < l, \quad 0 < y_2 < G(y_1) \}.
$

The purely periodic case can be addressed by somehow standard techniques in homogenization theory, as accomplished in \cite{R2, R3}. See \cite{R6,SP} for general references in homogenization and \cite{R7} for reticulated structures.    Observe that in this case the 
extension operators are very important for the convergence proof.

Moreover, if we assume that $g_\eps(\cdot)$ is independent of $\eps$, say $g_\eps(\cdot)=g(\cdot)$, and $h_0=\min_{x \in \R}\{h(x)\}$ then the variational formulation of the limit problem is:
\begin{equation}\label{HPC}
\begin{gathered}
\int^{1}_{0} \Big \{ \Big(g(x)+h_0 \Big) \; w_x(x)\, \varphi_x(x) + p(x)\, \omega(x)\, \varphi(x) \Big\} dx = \int^{1}_{0} \hat f(x) \, \varphi \, dx, \quad \forall \varphi\in H^1(0,1)
\end{gathered}
\end{equation}
where
$ p(x) = g(x) + \frac{1}{L_2} \int^{L_2}_{0}  h(s) \, ds, \hbox{  for all }x\in (0,1),$
 and the function $\hat f^\eps(x)=\int^{g(x)}_{-h_\eps(x)}f(x,y) \, dy$ satisfies that $\hat f^\eps\rightharpoonup \hat f$, w-$L^2(0,1)$.
 We refer to \cite{R4} for details. 
% In this case is not used any extension operator for the convergence proof,  is necessary to state a technical result about 
% the solutions of elliptics pde's  in certain rectangles that we  will show in Section 2.
% 
 In this work, we want to analyze the case where the thin domain is a region between two functions with different order of the oscillations.
 
 Our case is a combination of these two cases since both $g_\eps$ and $h_\eps$ are present. And we want to understand the effect of both terms at the same time in the limit equation.  Notice that the techniques used to solve each case separately are different so we will need to combine both techniques to get the limit problem in our case.  The main difference of the present work in relation to previous existing work in the literature, see for instance \cite{BFF, MP2, BZ} and references therein,  is that we allow two different order of oscillations in the boundary of the thin domain. 
 
 In Section 2 we state the notation and the problem that we will study. Furthermore, we are going to construct an extension operator that will be very important in the proof of the
convergence result. Finally, we state the main convergence result.  
 
 In Section 3 we rigorously prove the convergence result. In order to do so, we combine two different techniques: we use an extension operator 
 in the upper boundary and we define suitable rectangles in the lower boundary to apply the estimates that we obtained in Lemma \ref{basic-lemma} .

\section{Notation and statement of main result}
%Notation and important facts}
\label{basics}

To study the convergence of the solutions of (\ref{OPI0}) we first perform the change of variables $(x,y) \to (x, \eps y)$, which transforms the domain $R^\eps$ into the domain $\Omega^\eps$
\begin{equation} \label{domain}
\Omega^\epsilon = \Big\{ (x_1,x_2) \in \R^2 \; | \;  x_1 \in (0,1), \; -h(x_1/\eps^\alpha) < x_2 < g(x_1/\eps) \Big\}.
\end{equation}
Under this transformation, we obtain the equivalent linear elliptic problem 
\begin{equation} \label{P}
\left\{
\begin{gathered}
 - \frac{\partial^2 u^\epsilon}{{\partial x_1}^2} - 
\frac{1}{\epsilon^2} \frac{\partial^2 u^\epsilon}{{\partial x_2}^2} + u^\epsilon = f^\epsilon
\quad \textrm{ in } \Omega^\epsilon, \\
\frac{\partial u^\epsilon}{\partial x_1} \nu_1^\epsilon + \frac{1}{\epsilon^2} \frac{\partial u^\epsilon}{\partial x_2}\nu_2^\epsilon = 0
\quad \textrm{ on } \partial \Omega^\epsilon,
\end{gathered}
\right.
\end{equation}
where $f^\eps\in L^2(\Omega^\eps)$  satisfies 
$\| f^\epsilon \|_{L^2(\Omega^\epsilon)} \le C,$
for some $C > 0$ independent of $\epsilon$, and $\nu^\eps=(\nu_1^\eps,\nu_2^\eps)$ is the outward
unit normal to $\partial\Omega^\eps$.
Observe that $\Omega^\eps$ is not a thin domain anymore but there appears a  factor $1/\eps^2$ in front of the derivative in the $x_2$. Moreover, the domain  has very wild oscillatory behavior at the top and bottom boundary.

For the analysis we will construct an extension operator for functions defined in the set $\Omega^\eps$, but which will extend the function only over the upper part of the boundary. Hence, let us consider the following open set:
\begin{align} \label{domain}
{\widetilde{\Omega}}^\epsilon = \Big\{ (x_1,x_2) \in \R^2 \; | \;  x_1 \in (0,1), \; -h(x_1/\eps^\alpha) < x_2 < g_1 \Big\}.
\end{align}
\begin{lemma}\label{extension1}
With the notation above, there exists an extension operator
$$
P_\eps \, \in \, \mathcal{L} ( L^p(\Omega^\eps), \, L^p({\widetilde{\Omega}}^\epsilon)) \cap \mathcal{L} ( W^{1,p}(\Omega^\eps), \, W^{1,p}({\widetilde{\Omega}}^\epsilon))
$$
such that for any $\varphi \, \in \, W^{1,p}(\Omega^\eps)$,
\begin{eqnarray}\label{extension}
|| P_\eps \varphi ||_ { L^p({\widetilde{\Omega}}^\epsilon) } \le K ||  \varphi ||_ { L^p(\Omega^\eps)}, \; \left\| \frac{\partial P_\eps \varphi}{\partial x_2} \right\|_ { L^p({\widetilde{\Omega}}^\epsilon) } \le K  \left\| \frac{\partial \varphi}{\partial x_2} \right\| _ { L^p(\Omega^\eps)}\;\\
 \hbox{and} \; \;\left\| \frac{\partial P_\eps \varphi}{\partial x_1} \right\|_ { L^p({\widetilde{\Omega}}^\epsilon) } \le K \left \{ \left\| \frac{\partial \varphi}{\partial x_1} \right\| _ { L^p(\Omega^\eps)} + \eta(\eps) \left\| \frac{\partial \varphi}{\partial x_2} \right \| _ { L^p(\Omega^\eps)}\right\}\\
\end{eqnarray}
where $1 \le p \le \infty,$ $K$ a constant independent of $\eps$ and 
$\eta(\eps) = \sup_{x \in I}\{ |g'_\eps(x)|\}.  $
\end{lemma}
\emph{Proof.} The extension operator is constructed with a reflection procedure over the upper boundary, as in \cite{R2}. 
%Let us consider two cases. On the one hand, if we have that $g_1- h_0 \le 2 g_0$, which implies that $g(x_1/ \eps) - h(x_1/ \eps ^{\alpha}) \le 2g_0$, we can define the operator:
%$$
%(P_\eps \varphi)(x_1, x_2) = 
%\left\{
%\begin{array}{ll}
%\varphi(x_1, x_2) & (x_1, x_2) \in \Omega^\eps \\
%\varphi(x_1, 2g(x_1/ \eps) - x_2) & (x_1,x_2) \in \widetilde{\Omega}^\epsilon  / \Omega^\eps.
%\end{array}
%\right. 
%$$
%It easily follows that this operator satisfies the inequalities (\ref{extension}).
%On the other hand, if we are in the case where  $g_1- h_0 > 2 g_0$,  we will need to extend first the function $\varphi$ in the direction of negative $x_2$. Hence,
%we define a finite number of successive extensions such that, in the step $n$, we extend $\varphi$ to the open set  $\Omega_n^\epsilon = \{ (x_1,x_2) \in \R^2 \; | \;  x_1 \in (0,1), \; -h_0- ng_0 < x_2 < g(x_1/\eps) \}$ by reflecting
%$$
% \varphi_n(x_1, x_2) = 
%\left\{
%\begin{array}{ll}
%\varphi_{n-1}(x_1, x_2) & (x_1, x_2) \in \Omega_{n-1}^\eps \cup \Omega^\eps\\
%\varphi \big(x_1,  - x_2 -2\max \{h(x_1/\eps^\alpha), (n-1)g_0+h_0\}\big) & (x_1,x_2) \in (\Omega_n^\epsilon / \Omega_{n-1}^\epsilon)  / \big(\Omega^\eps \cap \Omega_n^\eps\big)
%\end{array}
%\right. .
%$$
%Choosing n large enough so that $ng_0 > g_1 -h_0$, constructing $\varphi_n$ and applying to $\varphi_n$ the procedure by reflection in the $x_2$ direction
%along the oscillating upper boundary, we obtain the extension operator which  satisfies the inequalities (\ref{extension}).
\qedsymbol

\par\medskip

Now, we state the convergence result:

\begin{theorem}\label{teorema}
Assume that $f^\eps \in L^2(\Omega^\eps)$  satisfies $\| f^\epsilon \|_{L^2(\Omega^\epsilon)} \le C$ with $C$ independent of the parameter $\eps$ and that there exists $\hat f\in L^2(0,1)$ such that  $
\hat f^\eps \rightharpoonup \hat f, \; w-L^2(0,1)$, where $\hat f^\eps (x_1) \equiv \int_{-h(x_1/ \eps^{\alpha})}^{g_1} \widetilde{f^\eps}(x_1,x_2)\, dx_2$.
 Let $u^\eps$ be the unique solution of (\ref{P}). Then, there exists $u_0 \in H^1(0,1)$ such that if $P_\eps$ is the extension operator constructed in Lemma \ref{extension1}, we have
 $ \|P_\eps u^\eps -u_0\|_{L^2(\widetilde\Omega^\eps)}\to 0$ and $u_0$   is the unique weak solution of the Neumann problem
 \begin{equation}\label{limit-equation-variational}
\int_{0}^{1} \Big\{\hat q \frac{\partial u_0}{\partial x_1}\frac{\partial \varphi}{\partial x_1} 
+ (\frac{|Y^*|}{L_1} +p)\, u_0 \,\varphi \,\Big\} dx_1 = \int_0^1  \hat f \, \varphi \, dx_1,\quad \forall \varphi\in H^1(0,1).
\end{equation}
where $Y^*$ is the basic cell
$$Y^*=\{ (y_1,y_2) \in \R^2 \; : \; 0<y_1<L_1 \textrm{ and } -h_0<y_2<g(y_1)\}.$$
The homogenized constant coefficients are defined by 
\begin{equation}\label{definition-q-hat}
\hat q \equiv \int_{-h_0}^{g_1} q(s) \,ds= \frac{1}{L_1}\int_{Y^*} \Big\{ 1- \frac{\partial X}{\partial y_1}(y_1, y_2) \Big\} dy_1 dy_2, 
\quad p= \frac{1}{L_2}\int_0^{L_2} h(s) ds - h_0,
\end{equation}
 where $X$ is the unique solution (up to constants) which is $L_1$-periodic in the first variable,  of the problem: 
\begin{equation} \label{AUX}
\left\{
\begin{gathered}
- \Delta X  =  0  \textrm{ in } Y^* \\
\frac{\partial X}{\partial N}  =  0  \textrm{ on } B_2  \\
\frac{\partial X}{\partial N} 
= - \frac{g'(y_1)}{\sqrt{1+g'(y_1)^2}} \textrm{ on } B_1
% \\
% X(0,y_2) = X(L_1,y_2) \textrm{ on } B_0\\
%   \int_{Y^*} X \; dy_1 dy_2 = 0  
\end{gathered}
\right.
\end{equation}
$B_0$ is the lateral part of the  boundary, $B_1$ is the upper boundary  and $B_2$ is the lower boundary of $\partial Y^*$.\\
\end{theorem}

\begin{remark}
If the non homogeneous term $f^\eps(x_1,x_2)$ is a fixed function depending only on the first variable, that is, $f^\eps(x_1,x_2)=f(x_1)$, it is easy to see that $\hat f(x_1)=(\frac{|Y^*|}{L_1} +p)f(x_1)$ and therefore, \eqref{limit-equation-variational} is the variational version of 
\begin{equation} \label{limit-equation1}
\left\{
\begin{gathered}
-\frac{\hat q}{\frac{|Y^*|}{L_1} +p}w_{xx} + w =  f(x), \quad x \in (0,1)\\
w'(0)=w'(1)=0
\end{gathered}
\right.
\end{equation} 
Notice that in case $h(\cdot)\equiv 0$, then $p=0$ and  $\frac{\hat q}{|Y^*|/L_1}=\frac{1}{|Y^*|}\int_{Y^*}(1-\frac{\partial X}{\partial y_1})=q_0$ and we recover \eqref{PC}.

\end{remark}

\section{Proof of the main result}

The variational formulation of (\ref{P}) is:  find $u^\epsilon \in H^1(\Omega^\epsilon)$ such that 
\begin{equation} \label{VFP}
\int_{\Omega^\epsilon} \Big\{ \frac{\partial u^\epsilon}{\partial x_1} \frac{\partial \varphi}{\partial x_1} 
+ \frac{1}{\epsilon^2} \frac{\partial u^\epsilon}{\partial x_2} \frac{\partial \varphi}{\partial x_2}
+ u^\epsilon \varphi \Big\} dx_1 dx_2 = \int_{\Omega^\epsilon} f^\epsilon \varphi dx_1 dx_2, 
\quad \forall \varphi \in H^1(\Omega^\epsilon).
\end{equation}
Taking $\varphi = u^\epsilon$ in  expression (\ref{VFP}) and using that $\|f^\eps\|_{L^2(\Omega^\eps)}\leq C$,  we easily obtain the a priori bounds
\begin{equation} \label{EST0}
\begin{gathered}
\| u^\epsilon \|_{L^2(\Omega^\epsilon)}, \Big\| \frac{\partial u^\epsilon}{\partial x_1} \Big\|_{L^2(\Omega^\epsilon)}
\textrm{ and } \frac{1}{\epsilon} \Big\| \frac{\partial u^\epsilon}{\partial x_2} \Big\|_{L^2(\Omega^\epsilon)} 
\le C.
\end{gathered}
\end{equation}

If we denote by $\; \widetilde{}\,\,$ the standard extension by zero and by $\chi^\eps$  the characteristic function  of $\Omega^\eps$, we may write (\ref{VFP}) as

\begin{equation}\label{VFP2}
\int_{\Omega_0} \Big\{\widetilde{ \frac{\partial u^\epsilon}{\partial x_1}} \frac{\partial \varphi}{\partial x_1} 
+ \frac{1}{\epsilon^2}\widetilde{ \frac{\partial u^\epsilon}{\partial x_2}} \frac{\partial \varphi}{\partial x_2} \Big\} 
+ \int_{\widetilde \Omega^\eps_-} \Big\{\widetilde{ \frac{\partial u^\epsilon}{\partial x_1}} \frac{\partial \varphi}{\partial x_1} 
+ \frac{1}{\epsilon^2} \widetilde{\frac{\partial u^\epsilon}{\partial x_2}} \frac{\partial \varphi}{\partial x_2} \Big\} 
+ \int_{\widetilde\Omega^\epsilon} \chi^\eps P_\eps u^\epsilon \varphi \,  = \int_{\widetilde\Omega^\epsilon} \chi^\eps f^\epsilon \varphi \,
\, \forall \varphi \in H^1(\Omega^\epsilon), 
\end{equation}
%\begin{eqnarray} \label{VFP2}
%& & \int_{\widetilde \Omega^\epsilon_+} \Big\{\widetilde{ \frac{\partial u^\epsilon}{\partial x_1}} \frac{\partial \varphi}{\partial x_1} 
%+ \frac{1}{\epsilon^2}\widetilde{ \frac{\partial u^\epsilon}{\partial x_2}} \frac{\partial \varphi}{\partial x_2} \Big\} dx_1 dx_2
%+ \int_{\widetilde \Omega^\eps_-} \Big\{\widetilde{ \frac{\partial u^\epsilon}{\partial x_1}} \frac{\partial \varphi}{\partial x_1} 
%+ \frac{1}{\epsilon^2} \widetilde{\frac{\partial u^\epsilon}{\partial x_2}} \frac{\partial \varphi}{\partial x_2} \Big\} dx_1 dx_2
%\nonumber \\
%& &  + \int_{\widetilde\Omega^\epsilon} \chi^\eps P_\eps u^\epsilon \varphi \, dx_1 dx_2 = \int_{\widetilde\Omega^\epsilon} \chi^\eps f^\epsilon \varphi \, dx_1 dx_2, 
%\, \forall \varphi \in H^1(\Omega^\epsilon).
%\end{eqnarray}
%
where we divide the domain $\widetilde\Omega^\eps$ in two parts: one of them,  $\widetilde\Omega^\eps_-$, carries all the oscillations and the other $\Omega_0$ is a fixed domain, that is, 
\begin{align} \label{DOMAINS}
\widetilde \Omega^\eps_-=  \{ (x_1, x_2) \in \R^2  | \, x_1 \in (0,1), \, -h(x_1/\eps^{\alpha}) < x_2 < -h_0 \} \nonumber \\
\Omega_0  =  \{ (x_1, x_2) \in \R^2  | \, x_1 \in (0,1), \, -h_0 < x_2 < g_1 \}. 
\end{align}

Before we start with the proof of the main result, let us state some relevant estimates on the solutions of certain elliptics problems, posed in rectangles of the type  
\begin{equation}\label{rectangles}
Q_\eps=\{ (x,y) \in \R^2 \; | \;  -\eps^\alpha<x<\eps^\alpha, \, 0<y<1\}, \; \hbox{with}\; \alpha>1.
\end{equation} 
As a matter of fact, 
for $u_0(\cdot)\in H^1(-\eps^\alpha,\eps^\alpha)$, we define the function $u^\eps(x,y)$  as the unique solution of
\begin{equation} \label{P-basic}
\left\{
\begin{gathered}
 - \frac{\partial^2 u^\epsilon}{{\partial x}^2} - 
\frac{1}{\epsilon^2} \frac{\partial^2 u^\epsilon}{{\partial y}^2}= 0
\quad \textrm{ in } Q_\epsilon, \\
\qquad u(x,0)=u_0(x),\quad  \textrm{ on } \Gamma_\eps,\\
\frac{\partial u}{\partial \nu}=0,\quad  \textrm{ on } \partial Q_\eps\setminus \Gamma_\eps
\end{gathered}
\right.
\end{equation}
where $\nu$ is the outward unit normal to $\partial Q_\eps$ and 
$
\Gamma_\eps = \{ (x,0) \in \R^2 \, | \, -\eps^\alpha<x<\eps^\alpha \}.
$

We have the following,
\begin{lemma}\label{basic-lemma}
With the notation from above, if we denote by $\bar u_0$ the average of $u_0$ in $\Gamma_\eps$, that is 
$
\bar u_0=\frac{1}{2\eps^\alpha}\int_{-\eps^\alpha}^{\eps^\alpha}u_0(x) \, dx
$
then there exists a constant $C$,  independent of $\eps$ and $u_0$, such that

\begin{equation}\label{estimate-L2}
\int_0^1\int_{-\eps^\alpha}^{\eps^\alpha}|u^\eps(x,y)-\bar u_0|^2 \, dxdy \leq C\eps^{\alpha-1} \|u_0\|_{L^2(-\eps^\alpha,\eps^\alpha)}^2
\end{equation}
and
\begin{equation}\label{basic-estimate}
\left\|\frac{\partial u^\eps}{\partial x}\right\|_{L^2(Q_\eps)}^2+\frac{1}{\eps^2}\left\|\frac{\partial u^\eps}{\partial y}\right\|^2_{L^2(Q_\eps)}\leq C \eps^{\alpha -1}
\left\|\frac{\partial u_0}{\partial x}\right\|_{L^2(-\eps^\alpha,\eps^\alpha)}^2.
\end{equation}
\end{lemma}
\emph{Proof.} See \cite{R4} for details. \qedsymbol

%We define now an appropriate partition for the unit interval $[0,1]$ which is related to the function $h_\eps$ and which will allow us to analyze in detail the effect of the oscillations at the bottom in the limit equation. 
%%use to
%%deadefine rectangles of the type (\ref{rectangles}).
%Let us denote by $N_\eps$ the largest integer such that $N_\eps L_2\eps^{\alpha}<1$, where $L_2$ is the period of the function $h$. 
%Observe that $N_\eps\sim L_2^{-1} \eps^{-\alpha}$. Let 
%\begin{equation} \label{eqG00}
%h_{n,\eps}=\min_{x\in [(n-1)L_2 \eps^{\alpha}, nL_2 \eps^{\alpha}]} h\left(\frac{x}{\eps^\alpha}\right), \quad n=1,2\ldots, N_\eps
%\end{equation}
%and $\gamma_{n,\eps}\in [(n-1)L_2 \eps^{\alpha}, nL_2 \eps^{\alpha}]$ a point where the minimum (\ref{eqG00}) is attained, 
%that is, $h(\frac{ \gamma_{n,\eps}}{\eps^\alpha})=h_{n,\eps}$ where $\gamma_{n,\eps}$ does not need to be uniquely defined.
%By extension, let us denote by $\gamma_{0,\eps}=0$ and $\gamma_{ N_\eps+1,\eps}=1$.
%
%Note that the set
%$
%\{ \gamma_{0,\eps}, \gamma_{1,\eps}, ..., \gamma_{N_\eps + 1,\eps} \}
%$
%defines a partition for the unit interval $[0,1]$. Moreover, due to that $h(\cdot)$ is $L_2-$periodic we have that $h_{n,\eps} = h_0 \; \hbox{for} \; n=1,2\ldots, N_\eps.$
%\section{Proof of the main result}
%In this section we obtain the limit of the sequence $\{ u^\eps\}_{\eps>0}$ given by the Neumann problem (\ref{P}) as the parameter $\eps$ goes to zero. 

\par\bigskip
\noindent\emph{Proof of Theorem \ref{teorema}.}
The idea is to pass to the limit in (\ref{VFP2}) constructing appropriate test functions. 
First, we study the limit of the different functions that form the integrands of (\ref{VFP2}).
\par\medskip\noindent {\bf  (a). Limit in the extended functions. }
Using the a priori estimate (\ref{EST0}) and the results from Lemma \ref{extension1} we obtain that $P_\eps u^\eps|_{\Omega_0} \in H^1(\Omega_0)$ and we can extract a subsequence of $\{ P_\eps u^\eps|_{\Omega_0} \} \subset H^1(\Omega_0)$, denoted again by $P_\eps u^\epsilon$, such that
\begin{equation} \label{WC0}
\begin{gathered}
P_\eps u^\epsilon \rightharpoonup u_0 \quad w-H^1(\Omega_0) \\
%P_\eps u^\epsilon \rightarrow u_0 \quad s-H^s(\Omega_0) \textrm{ for all } s \in [0, 1) \textrm{ and }\\
\frac{\partial P_\eps u^\epsilon}{\partial x_2} \rightarrow 0 \quad s-L^2(\Omega_0) \\
\end{gathered}
\end{equation}
as $\epsilon \to 0$ for some  $u_0 \in H^1(\Omega_0)$.

A consequence of the limits (\ref{WC0}) is that $ u_0(x_1,x_2)$ does not depend on the variable $x_2$. 
Moreover, we have that the restriction of $P_\eps u^\eps$ to the coordinate axis $x_1$ converges to $u_0$.
That is, 
$
P_\eps u^\eps|_{\Gamma} \rightarrow u_0 \quad s-H^s(\Gamma)
$
for all $s \in [0,1/2)$ where  $\Gamma = \{ (x_1,0) \in \R^2 \, | \, x_1 \in (0,1) \}$. 
Consequently, we obtain
$
\| P_\eps u^\eps - u_0 \|_{L^2(\Gamma)} \to 0 \textrm{ as } \eps \to 0.
$\\
 In view of the above limit, one has the $L^2-$convergence of $P_\eps u^\eps $ to $u_0$, that is
\begin{equation} \label{L2CONV}
\| P_\eps u^\eps - u_0 \|_{L^2(\widetilde\Omega^\eps)} \to 0 \textrm{ as } \eps \to 0.
\end{equation}
In fact, on the one hand we have 
\begin{eqnarray*}
\| P_\eps u^\eps(x_1,0) - u_0(x_1) \|_{L^2(\widetilde \Omega^\eps)}^2  =
\int_0^1 \int_{-h(x_1/\eps^{\alpha})}^{g_1} | P_\eps u^\eps(x_1,0) - u_0(x_1) |^2 \, dx_2 dx_1\\
\leq C(g,h) \, \| P_\eps u^\eps - u_0 \|_{L^2(\Gamma)} \to  0 \textrm{ as } \eps \to 0.
\end{eqnarray*}
On the other hand,
\begin{eqnarray*}
 \| P_\eps u^\eps(x_1,x_2) - P_\eps u^\eps(x_1,0) \|_{L^2(\widetilde\Omega^\eps)}^2 = 
\int_0^1 \int_{-h(x_1/\eps^{\alpha})}^{g_1} | P_\eps u^\eps(x_1,x_2) - P_\eps u^\eps(x_1,0) |^2 \, dx_1 dx_2 \\
 \leq  \int_0^1 \int_{-h(x_1/\eps^{\alpha})}^{g_1} \left(  \int_0^{x_2} \left| \frac{\partial P_\eps u^\eps}{\partial x_2}(x_1,s) \right|^2 ds \right) 
\, |x_2| \, dx_2 dx_1
\leq  C(h,g) \, \left\| \frac{\partial P_\eps u^\eps}{\partial x_2} \right\|_{L^2(\widetilde\Omega^\eps)}^2 \leq  \eps \, \hat C(h,g)  \to 0 \textrm{ as } \eps \to 0.
\end{eqnarray*}
Finally
\begin{eqnarray*}
\| P_\eps u^\eps - u_0 \|_{L^2(\widetilde\Omega^\eps)} & \leq & \| P_\eps u^\eps(x_1,x_2) - P_\eps u^\eps(x_1,0) \|_{L^2(\widetilde\Omega^\eps)} + 
\|P_\eps u^\eps(x_1,0) - u_0(x_1) \|_{L^2(\widetilde\Omega^\eps)} \to  0,\\
 \textrm{ as } \eps \to 0.
\end{eqnarray*}

\par\noindent {\bf  (b). Limit in the tilde functions. }

 From the a priori estimates (\ref{EST0}) we know that there exists a function $\xi^* \in L^2(\Omega_0)$, such that, up to subsequences
\begin{equation} \label{WC1}
\begin{gathered}
\widetilde{\frac{\partial u^\epsilon}{\partial x_1}} \rightharpoonup \xi^* \; w-L^2(\Omega_0) \;\; \hbox{and} \;\;
\widetilde{\frac{\partial u^\epsilon}{\partial x_2}} \rightarrow 0 \; s-L^2(\Omega_0); \quad \hbox{as}\; \epsilon \to 0.
\end{gathered}
\end{equation}

\par\noindent {\bf  (c). Limit of $\chi^\eps$. }

Let $\chi$ be the characteristic function of the representative cell $Y^*$.
We extend $\chi$ periodically on the variable $y_1 \in \R$ and denote this extension again by $\chi$.
Clearly, by construction,
$
\chi^\epsilon(x_1,x_2) 
= \chi (x_1/\epsilon,x_2),\hbox{ for }(x_1,x_2)  \in \Omega^\epsilon_+.
$\\
Consequently, by the Average Theorem and the Lebesgue's Dominated Convergence Theorem we obtain 

\begin{equation} \label{CHIAIM0}
\chi^\epsilon \stackrel{\eps\to 0}{ \rightharpoonup} \theta \quad w^*-L^\infty(\Omega_0), \quad \hbox{where} \;\;
\theta(x_2) := \frac{1}{L_1} \int_0^{L_1} \chi(s,x_2) ds  \quad\forall x_2\in (-h_0,g_1).
\end{equation}

%{\color{red}  
%\par\noindent {\bf  (d). Limit of $f^\eps$. }
%
%From hypotheses $\|f^{ \eps}\|_{L^2(\Omega^{ \eps})}\leq C$ independent of $\eps$. Hence, $\|\widetilde{f^{ \eps}}\|_{L^2(\widetilde\Omega^{ \eps})}\leq C$
%and the function $\hat f^\eps$ defined by 
%\begin{equation} \label{FHF}
%\hat f^\eps (x_1) \equiv \int_{-h(x_1/ \eps^{\alpha})}^{g_1} \widetilde{f^\eps}(x_1,x_2)dx_2
%\end{equation}
%belongs to $L^2(0,1)$ and satisfies 
%$
%\| \hat f^\eps \|_{L^2(0,1)} \leq C(h,g) 
%$ 
%for some constant $C(h,g)$ independent of $\eps$. Then, via subsequences, we have the existence of a function 
%$\hat f=\hat f(x_1)\in L^2(0,1)$ such that 
%\begin{equation}\label{LIMITF}
%\hat f^\eps \rightharpoonup \hat f \qquad w-L^2(0,1).
%\end{equation}
%
%}

\par\noindent {\bf(d) Test  functions. }

In order to construct appropriate test functions that will allow us to pass the limit in the variational formulation (\ref{VFP2}),  we are going to need to define a partition of the unit interval $[0,1]$ which is related to the function $h_\eps$ and which will allow us to analyze in detail the effect of the oscillations at the bottom in the limit equation. 
%use to
%deadefine rectangles of the type (\ref{rectangles}).
Hence, denote by $N_\eps$ the largest integer such that $N_\eps L_2\eps^{\alpha}<1$, where $L_2$ is the period of the function $h$. 
Observe that $N_\eps\sim L_2^{-1} \eps^{-\alpha}$. Let 
\begin{equation} \label{eqG00}
h_{n,\eps}=\min_{x\in [(n-1)L_2 \eps^{\alpha}, nL_2 \eps^{\alpha}]} h\left(\frac{x}{\eps^\alpha}\right), \quad n=1,2\ldots, N_\eps
\end{equation}
and $\gamma_{n,\eps}\in [(n-1)L_2 \eps^{\alpha}, nL_2 \eps^{\alpha}]$ a point where the minimum (\ref{eqG00}) is attained, 
that is, $h(\frac{ \gamma_{n,\eps}}{\eps^\alpha})=h_{n,\eps}$ where $\gamma_{n,\eps}$ does not need to be uniquely defined.
By extension, let us denote by $\gamma_{0,\eps}=0$ and $\gamma_{ N_\eps+1,\eps}=1$.

Note that the set
$
\{ \gamma_{0,\eps}, \gamma_{1,\eps}, ..., \gamma_{N_\eps + 1,\eps} \}
$
defines a partition for the unit interval $[0,1]$. Moreover, due to that $h(\cdot)$ is $L_2-$periodic we have that $h_{n,\eps} = h_0 \; \hbox{for} \; n=1,2\ldots, N_\eps.$

We define now  the test functions as follows. 
With $\phi \in H^1(0,1)$, we consider $\varphi^\eps\in H^1(\widetilde \Omega^\eps)$ defined as
\begin{equation} \label{TESTF}
\begin{gathered}
\varphi^\eps(x_1,x_2) = \left\{
\begin{array}{ll}
X^\eps_n(x_1,x_2), & (x_1,x_2) \in \widetilde\Omega^\eps_-\cap Q^\eps_n,\quad n=1,2,\ldots \\
\phi(x_1), & (x_1,x_2) \in \widetilde \Omega_+^\eps \equiv \Omega_0 
\end{array}
\right.
\end{gathered}
\end{equation}
where  $Q^\eps_n$ is the rectangle 
$
Q^\eps_n=\{ (x_1,x_2) \, | \, \gamma_{n,\eps}<x_1<\gamma_{n+1,\eps}, \, -h_1<x_2<-h_0 \}
$
and the function $X^\eps_n$ is the solution of the problem
\begin{equation} \label{AUXSOL}
\left\{
\begin{gathered}
- \frac{\partial^2 X_n^\eps}{\partial x_1^2} - \frac{1}{\eps^2} \frac{\partial^2 X_n^\eps}{\partial x_2^2} 
= 0, \quad \textrm{ in } Q^\eps_n \\
\frac{\partial X_n^\eps}{\partial N^\eps}=0, \quad \textrm{ on }  \partial Q^\eps_n \backslash \Gamma_n^\eps  \\
X_n^\eps(x_1,x_2) = \phi(x_1),  \quad \textrm{ on } \Gamma_n^\eps
\end{gathered}
\right. 
\end{equation}
where $\Gamma_n^\eps$ is the base of the rectangle, that is, 
$
\Gamma_n^\eps = \{ (x_1, -h_0): \gamma_{n,\eps}\leq x_1\leq \gamma_{n+1,\eps}\}.
$

From Lemma \ref{basic-lemma} we have
\begin{equation} \label{ESTX}
\left\| \frac{\partial X_n^\eps}{\partial x_1^2} \right\|^2_{L^2(Q^\eps_n)}
+ \frac{1}{\eps^2} \left\| \frac{\partial X_n^\eps}{\partial x_2^2} \right\|^2_{L^2(Q^\eps_n)}
\leq C \eps^{\alpha-1} \| \phi' \|^2_{L^2(\gamma_{n,\eps},\gamma_{n+1,\eps})}.
\end{equation}

Furthermore, since
$$
\varphi^\eps(x_1,x_2) - \phi(x_1) = \varphi^\eps(x_1,x_2) - \varphi^\eps(x_1,0) 
= \int_0^{x_2} \frac{\partial \varphi^\eps}{\partial x_2}(x_1,s) \, ds,
$$
we have by (\ref{TESTF}) and (\ref{ESTX}) that

\begin{equation} \label{TFCONV}
\| \varphi^\eps - \phi \|_{L^2(\widetilde\Omega^\eps)} \to 0 \textrm{ as } \eps \to 0.
\end{equation}

\par\noindent {\bf  (e) Passing to the limit. }

We can now pass to the limit in (\ref{VFP2}) by making use of test functions $\varphi^\eps$ defined above. For this, we study the convergence of each term in (\ref{VFP2}).
\begin{itemize}

\item First integrand:
\begin{equation} \label{INT2}
\int_{\Omega_0} \Big\{ \widetilde{\frac{\partial u^\epsilon}{\partial x_1}} \frac{\partial \varphi^\eps}{\partial x_1} 
+ \frac{1}{\epsilon^2} \widetilde{\frac{\partial u^\epsilon}{\partial x_2}} \frac{\partial \varphi^\eps}{\partial x_2} \Big\} dx_1 dx_2
\to \int_{\Omega_0} \xi^*(x_1,x_2)  \phi'(x_1) \, dx_1 dx_2 \textrm{ as } \eps \to 0.
\end{equation}

Thanks to the choice of the test function (\ref{TESTF}) and the convergence (\ref{WC1}), we easily get  (\ref{INT2}).
\par
\item Second integrand: 
\begin{equation} \label{INT1}
\int_{\widetilde\Omega^\epsilon_-} \Big\{ \widetilde{ \frac{\partial u^\epsilon}{\partial x_1}} \frac{\partial \varphi^\eps}{\partial x_1} 
+ \frac{1}{\epsilon^2} \widetilde {\frac{\partial u^\epsilon}{\partial x_2}} \frac{\partial \varphi^\eps}{\partial x_2} \Big\} dx_1 dx_2
\to 0 \textrm{ as } \eps \to 0.
\end{equation}
From the definition of $\varphi^\eps$, the Cauchy-Schwarz inequality and the inequality (\ref{ESTX}) we have (\ref{INT1}).
\par
\item Third integrand: 
\begin{equation} \label{INT3}
\int_{\widetilde\Omega^\eps}\chi^\eps P_\eps u^\epsilon  \, \varphi^\eps \, dx_1 dx_2 \to \int_0^1 p \, u_0(x_1) \, \phi(x_1) \, dx_1+ \int_{\Omega_0} \theta(x_2) \, u_0(x_1) \, \phi(x_1) \, dx_1 dx_2  \textrm{ as } \eps \to 0
\end{equation}
where the constant $p$ is given by
$p= \frac{1}{L_2}\int_0^{L_2} h(s) ds - h_0.$

For this, note that we can rewrite the integral of the left side of (\ref{INT3}) as
\begin{eqnarray*}
\int_{\widetilde\Omega^\eps}\chi^\eps P_\eps u^\epsilon  \, \varphi^\eps \, dx_1 dx_2  =  
\int_{\widetilde\Omega^\eps}\chi^\eps \left( P_\eps u^\eps - u_0 \right) \, \varphi^\eps \, dx_1 dx_2\\
+ \int_{\widetilde\Omega^\eps} \chi^\eps u_0 \, \left( \varphi^\eps - \phi \right) \, dx_1 dx_2
+ \int_{\widetilde\Omega^\eps}  \chi^\eps  u_0 \, \phi \, dx_1 dx_2.
\end{eqnarray*}
From (\ref{L2CONV}) and (\ref{TFCONV}), we have that the first two terms in the right hand side above go to 0. 
%$$
%\begin{gathered}
%\int_{\widetilde\Omega^\eps}\chi^\eps \left( P_\eps u^\eps - u_0 \right) \, \varphi^\eps \, dx_1 dx_2 \to 0 \quad \textrm{ and } \int_{\widetilde\Omega^\eps} \chi^\eps u_0 \, \left( \varphi^\eps - \phi \right) \, dx_1 dx_2 \to 0, \hbox{ as } \eps \to 0
%\end{gathered}
%$$
Moreover, since
\begin{eqnarray*}
\int_{\widetilde\Omega^\eps}\chi^\eps  u_0  \, \phi \, dx_1 dx_2 = \int_{\widetilde\Omega^\eps_-} u_0  \, \phi \, dx_1 dx_2 + \int_{\Omega_0} \chi^\eps u_0  \, \phi \, dx_1 dx_2\\
= \int_0^1 u_0  \, \phi\left( h\left(\frac{x_1}{\eps^{\alpha}}\right) - h_0 \right) \, dx_1 + \int_{\Omega_0} \chi^\eps u_0  \, \phi \, dx_1 dx_2
\end{eqnarray*}
we get (\ref{INT3}) from the Average Theorem and  (\ref{CHIAIM0}) .

\item Fourth integrand: 
\begin{equation} \label{INT4}
\int_{\widetilde\Omega^\eps} \widetilde{f^\eps} \, \varphi^\eps \, dx_1 dx_2 \to \int_0^1  \hat f(x_1) \, \phi(x_1) \, dx_1 \textrm{ as } \eps \to 0.
\end{equation}

From (\ref{TFCONV}) and the hypotheses of the theorem we have (\ref{INT4}).
\end{itemize}

Therefore, using the convergences (\ref{INT1}), (\ref{INT2}), (\ref{INT3}) and (\ref{INT4}), we obtain the following limit variational formulation: 
\begin{eqnarray} \label{limitP}
\int_{\Omega_0} \left\{ \xi^*(x_1,x_2) \, \phi'(x_1) + \theta(x_2) \, u_0(x_1) \, \phi(x_1)\right\} dx_2 dx_1 + \int_0^1 p \, u_0(x_1) \, \phi(x_1) \, dx_1\nonumber\\
= \int_0^1 \hat f(x_1) \, \phi(x_1) \, dx_1, \, \forall \phi \in H^1(0,1).
\end{eqnarray}
with
$p= \frac{1}{L_2}\int_0^{L_2} h(s) ds - h_0.$

At this point the question is how to relate $u_0$ to $\xi^*$. In the following subsection we will show a equation for $\xi^*$.
\par\medskip\noindent {\bf  (f) Relation between $\xi^*$ and $u_0$. }

Let us consider the following families of isomorphisms 
$T^{\epsilon}_k : A^{\epsilon}_k \mapsto Y$ given by
\begin{equation} \label{ISO}
T^{\epsilon}_k(x_1,x_2) = (\frac{x_1 - \epsilon k L_1}{\epsilon},x_2)
\end{equation}
where 
$$
\begin{gathered}
A^{\epsilon}_k = \{ (x_1,x_2) \in \R^2 \; | \;  
\epsilon  k L_1 \leq x_1 < \epsilon L_1 (k+1), -h_0< x_2 < g_1 \}  \quad \hbox{and} \quad Y = (0, L_1) \times (-h_0, g_1).
\end{gathered}
$$
with $k \in \N$.
We can considerer extension operators $
P \in \mathcal{L}(H^1(Y^*),H^1(Y)) \cap \mathcal{L}(L^2(Y^*),L^2(Y))$,the proof is done in \cite{R4}.
Using these operators, the isomorphism (\ref{ISO}) and the unique solution of the auxiliary problem (\ref{AUX}) we define $\omega^{\epsilon}_k$ in $(x_1, x_2) \in A^{\epsilon}_k$ by
\begin{eqnarray*}
\omega^{\epsilon}_k (x_1,x_2) & = & x_1 - \epsilon \Big(  PX \circ T^\epsilon_k (x_1,x_2) \Big)
 =  x_1 - \epsilon \Big( P X (\frac{x_1 - \epsilon L_1 k}{\epsilon},x_2) \Big).
\end{eqnarray*}
Observe that for any $(x_1, x_2) \in \widetilde\Omega^{\eps}_+$ there is $k$ such that $(x_1, x_2) \in A^{\epsilon}_k$. 
Therefore, the function $\omega^{\epsilon}(x_1,x_2)=\omega^{\epsilon}_k (x_1,x_2)$ is well defined and $\omega^{\epsilon} \in H^1(\widetilde\Omega^{\eps}_+)$.
We introduce now the vector $\eta^\epsilon = (\eta_1^\epsilon,\eta_2^\epsilon)$ defined by 
\begin{equation} \label{VETA}
\eta_i^\epsilon(x_1,x_2) = \frac{\partial \omega^\epsilon}{\partial x_i}(x_1,x_2),  \quad (x_1,x_2)\in \Omega^\eps_+
\end{equation}
where $\Omega^\eps_+=\{ (x_1,x_2) \in \R^2 \; : \; 0<x_1<1 \textrm{ and } -h_0<x_2<g(x_1/\eps)\}.$

Taking into account the definition of $X$  if we consider a  test function $\psi \in H^1(\Omega^\eps_+)$ with $\psi=0$ in  neighborhood of the lateral boundaries, we get
\begin{equation}\label{equation-final}
\int_{\Omega^\eps_+}  \left(\eta_1^\epsilon   \frac{\partial \psi}{\partial x_1}  
+  \eta_2^\epsilon \dfrac{1}{\epsilon^2} \frac{\partial \psi}{\partial x_2} \right) dx_1 dx_2=0.
\end{equation}
Then,  with the variational formulation (\ref{VFP}) and the identity (\ref{equation-final}) we can write:
\begin{eqnarray} \label{VFPE2}
 \int_{\overline{\Omega}^\epsilon} \Big\{ \widetilde{\frac{\partial u^\epsilon}{\partial x_1}} \frac{\partial \varphi}{\partial x_1} 
+ \frac{1}{\epsilon^2}\widetilde{ \frac{\partial u^\epsilon}{\partial x_2}} \frac{\partial \varphi}{\partial x_2}
+ \chi^\eps P_\eps u^\epsilon \varphi \Big\} dx_1 dx_2 - \int_{\Omega^\eps_+}  \left(\eta_1^\epsilon   \frac{\partial \psi}{\partial x_1}  
+  \eta_2^\epsilon \dfrac{1}{\epsilon^2} \frac{\partial \psi}{\partial x_2} \right) dx_1 dx_2\nonumber\\
 = \int_{\widetilde{\Omega}^\epsilon} \chi^\eps f^\epsilon \varphi dx_1 dx_2, 
\; \forall \varphi \in H^1(\Omega^\epsilon).
\end{eqnarray}

We would like to pass to the limit in this expression. For this, we will construct appropriate test functions, which used in the identity (\ref{VFPE2}) allow us to pass to the limit in all
the terms.
\par\noindent {\bf  (g) Limit of $\omega^\eps$ and $\eta_1^\eps$.}
From the definition of $\omega^\eps$, we have  
\begin{equation} \label{OMEGAL}
\omega^\epsilon \to x_1, \; s-L^2(\Omega_0); \qquad 
\frac{\partial \omega^\epsilon}{\partial x_2} \to 0, \; s-L^2(\Omega_0);\qquad \;
\widetilde{\eta}_1^\epsilon \rightharpoonup  q,  \; w^*-L^\infty(\Omega_0), 
\end{equation}
where $$ q(x_2):= \frac{1}{L_1} \int_0^{L_1} 
\Big( 1 - \widetilde{\frac{\partial X}{\partial y_1}}  (s,x_2) \Big) \chi(s,x_2)ds $$
See \cite{R3} for more details.

\par\noindent {\bf  (h) Function test. }

Let $\phi= \phi(x_1) \in \mathcal{C}^\infty_0(0,1)$. 
%and consider the test function  $\phi(x_1) \omega^\epsilon (x_1, -h_0) \in H^1(0,1).$ 
We introduce the test function
\begin{equation} \label{TESTF1}
\begin{gathered}
\psi^\eps(x_1,x_2) = \left\{
\begin{array}{ll}
X^\eps_n(x_1,x_2), & (x_1,x_2) \in \widetilde\Omega^\eps_-\cap Q^\eps_n,\quad n=1,2,\ldots \\
\phi(x_1)\omega^\epsilon (x_1, x_2), & (x_1,x_2) \in \widetilde \Omega_+^\eps \equiv \Omega_0,  
\end{array}
\right.
\end{gathered}
\end{equation}
where $\omega_\eps$ is defined above and, as in (\ref{TESTF}), 
$Q^\eps_n$ is the rectangle $
Q^\eps_n=\{ (x_1,x_2) \, | \, \gamma_{n,\eps}<x_1<\gamma_{n+1,\eps}, \, -h_1<x_2<-h_0 \}
$
and the function $X^\eps_n$ is the solution of the problem
\begin{equation} \label{AUXSOL}
\left\{
\begin{gathered}
- \frac{\partial^2 X_n^\eps}{\partial x_1^2} - \frac{1}{\eps^2} \frac{\partial^2 X_n^\eps}{\partial x_2^2} 
= 0, \quad \textrm{ in } Q^\eps_n \\
\frac{\partial X_n^\eps}{\partial N^\eps}=0, \quad \textrm{ on }  \partial Q^\eps_n \backslash \Gamma_n^\eps  \\
X_n^\eps(x_1,x_2) = \phi(x_1)\omega^\epsilon (x_1, -h_0),  \quad \textrm{ on } \Gamma_n^\eps
\end{gathered}
\right. 
\end{equation}
where $\Gamma_n^\eps$ is the base of the rectangle, that is, 
$
\Gamma_n^\eps = \{ (x_1, -h_0): \gamma_{n,\eps}\leq x_1\leq \gamma_{n+1,\eps}\}.
$

Moreover, we define the function $
X^\eps(x_1,x_2) = X^\eps_n(x_1,x_2) \textrm{ as } (x_1,x_2) \in Q^\eps_n \cap \overline\Omega^\eps_-.$ 

From Lemma \ref{basic-lemma} and using the  the properties of $ \omega^\eps$ we have that the function $X^\eps$ is $ H^1(\widetilde\Omega^\eps_-)$ and satisfies the following estimate
\begin{equation}\label{ESTX3}
\begin{array}{l}
\displaystyle \left\| \frac{\partial X^\eps}{\partial x_1} \right\|^2_{L^2(\widetilde \Omega^\eps_-)}
+ \frac{1}{\eps^2} \left\| \frac{\partial X^\eps}{\partial x_2} \right\|^2_{L^2(\widetilde\Omega^\eps_-)}
 \leq  C \, \eps^{\alpha -1} \left\| \frac{\partial\big(\phi(x_1)\omega^\epsilon (x_1, -h_0)\big)}{\partial{x_1}} \right\|^2_{L^2(0,1)} \leq  \widetilde{C} \, \eps^{\alpha -1}.
\end{array}
\end{equation}
where $\widetilde{C}$ denotes a constant independent of $\eps$.
 Now,  we can argue as in (\ref{TFCONV}) and we obtain
\begin{equation} \label{TFCONV1}
\| \psi^\eps - \phi \overline P \omega^\eps  \|_{L^2(\widetilde\Omega^\eps)} \to 0 \textrm{ as } \eps \to 0.
\end{equation} 
where $\overline P \omega^\eps$ is the function defined on the set  $\{ (x_1,x_2) \in \R^2 \; | \;  x_1 \in (0,1), \; -h_1< x_2 <g _1 \}$ using a extension operator
obtained by reflection in the negative vertical direction along the line $x_2=-h_0$. .Indeed, since
\begin{eqnarray*}
&\psi^\eps(x_1,x_2)- \phi(x_1)\overline P\omega^\eps(x_1, x_2)=\psi^\eps(x_1,x_2)- \phi(x_1)\omega^\eps(x_1, -x_2-2h_0)\nonumber \\
&\qquad \quad \quad =\psi^\eps(x_1,x_2) - \psi^\eps(x_1, -x_2-2h_0)= \int_{-x_2-2h_0}^{x_2}\frac{\partial \psi^\eps}{\partial x_2} (x_1,s) d_s \quad \text{for} \, (x_1, x_2) \in \widetilde\Omega^\eps_-
\end{eqnarray*}
we have by (\ref{TESTF1})
$$
\| \psi^\eps - \phi \overline P \omega^\eps  \|_{L^2(\widetilde\Omega^\eps)} \leq C(g,h)\, \left\| \frac{\partial \psi^\eps}{\partial x_2} \right\|_{L^2(\widetilde\Omega^\eps)}^2
= C(g,h)\, \left\| \frac{\partial \omega^\eps}{\partial x_2}\phi \right\|_{L^2(\widetilde\Omega^\eps_+)}^2 + C(g,h)\, \left\| \frac{\partial X^\eps}{\partial x_2} \right\|_{L^2(\widetilde\Omega^\eps_-)}^2\to 0 \; \hbox{as} \; \eps \to 0.
$$

\par\noindent {\bf  (i) Passing to the limit. }

Now we pass to the limit in the equality (\ref{VFPE2}) considering the test functions $\varphi= \psi^\eps$ and $\psi = \phi u^\eps$. 
\begin{itemize}
\item First integrand: 
\begin{equation} \label{INTE1}
\int_{\widetilde\Omega^\epsilon_-} \Big\{ \widetilde{ \frac{\partial u^\epsilon}{\partial x_1}} \frac{\partial \psi^\eps}{\partial x_1} 
+ \frac{1}{\epsilon^2} \widetilde {\frac{\partial u^\epsilon}{\partial x_2}} \frac{\partial \psi^\eps}{\partial x_2} \Big\} dx_1 dx_2
\to 0 \textrm{ as } \eps \to 0.
\end{equation}

Taking account the definition of $\psi^\eps$, the Cauchy-Schwarz inequality and the estimate (\ref{ESTX3}) we obtain the convergence (\ref{INTE1}).

\par 
\item Second integrand:
\begin{eqnarray}\label{INTE2}
  \int_{\widetilde{\Omega}^\epsilon_+} \Big\{ \widetilde{\frac{\partial u^\epsilon}{\partial x_1}} \frac{\partial \psi^\eps}{\partial x_1} 
+ \frac{1}{\epsilon^2}\widetilde{ \frac{\partial u^\epsilon}{\partial x_2}} \frac{\partial  \psi^\eps}{\partial x_2}
  - \eta_1^\epsilon   \frac{\partial (\phi u^\eps)}{\partial x_1}  
-  \eta_2^\epsilon \dfrac{1}{\epsilon^2} \frac{\partial (\phi u^\eps)}{\partial x_2}\Big\} dx_1 dx_2\nonumber \\
 \to \int_{\Omega_0} \Big\{ \xi^* \frac{\partial \phi}{\partial x_1} \, x_1 -    q\frac{\partial \phi}{\partial x_1} u_0
\Big\}  dx_1 dx_2  \textrm{ as } \eps \to 0.
\end{eqnarray}
 From the definitions of $\eta_i^\epsilon$ and $\psi^\eps$ the second integrand reduces to
$
\int_{\widetilde \Omega^\eps_+} \Big\{\widetilde{{\frac{\partial u^\epsilon}{\partial x_1}}} \frac{\partial \phi}{\partial x_1} \, \omega^\epsilon 
-  \widetilde{\eta_1^\epsilon}  \frac{\partial \phi}{\partial x_1} \,  P_\eps u^\epsilon 
\Big\}  dx_1 dx_2.
$

Therefore, using convergences (\ref{L2CONV}), (\ref{WC1}) and  (\ref{OMEGAL}), we have (\ref{INTE2}).
\par
\item Third integrand:
\begin{eqnarray} \label{INTE3}
\int_{\widetilde\Omega^\eps}\chi^\eps P_\eps u^\epsilon  \, \psi^\eps \, dx_1 dx_2 \to \int_0^1 p \, u_0(x_1) \, \phi(x_1) \, x_1 \,dx_1
+ \int_{\Omega_0} \theta(x_2) \, u_0(x_1) \, \phi(x_1) \, x_1\, dx_1 dx_2, \,  \textrm{ as } \eps \to 0.
\end{eqnarray}
Following along the lines of the proof of the convergence (\ref{INT3}) we have this convergence.

\par
\item Fourth integrand:
\begin{equation} \label{INTE4}
\int_{\widetilde\Omega^\eps} \widetilde{f^\eps} \, \psi^\eps \, dx_1 dx_2 \to \int_0^1  \hat f(x_1) \, \phi(x_1) \, x_1 dx_1 \textrm{ as } \eps \to 0.
\end{equation}
Using the same computations as those made to derive (\ref{INT4}) we obtain (\ref{INTE4})
\end{itemize}

 Now, by the convergences shown in (\ref{INTE1}), (\ref{INTE2}), (\ref{INTE3}) and (\ref{INTE4}), we can pass to the limit in (\ref{VFPE2}) 
 considering the test functions $\varphi= \psi^\eps$ and $\psi = \phi u^\eps$. More precisely, we have
\begin{align} \label{limitP1}
 \int_{\Omega_0} \Big\{ \xi^* \frac{\partial \phi}{\partial x_1} \, x_1 -    q\frac{\partial \phi}{\partial x_1} u_0
\Big\}  dx_1 dx_2   + \int_0^1 p \, u_0 \, \phi\, x_1 \, dx_1
+ \int_{\Omega_0} \theta \, u_0 \, \phi\, x_1\,dx_1 dx_2
= \int_0^1 \hat f \, \phi \, x_1 \,dx_1 \quad \forall \phi \in \mathcal{C}^\infty_0(0,1)
\end{align}
where  $p$ and $q$ are given by
$$p= \frac{1}{L_2}\int_0^{L_2} h(s) ds - h_0,\quad
q(x_2)=\frac{1}{L_1} \int_0^{L_1} 
\Big( 1 - \widetilde{\frac{\partial X}{\partial y_1}}  (s,x_2) \Big) \chi(s,x_2)ds.
$$

Taking the test function $\phi x_1$ in (\ref{limitP}) we obtain 
\begin{eqnarray} \label{limitP2}
\int_{\Omega_0} \xi^*(x_1,x_2) \, \frac{\partial}{\partial x_1}(\phi x_1) dx_2 dx_1 + \int_0^1 p \, u_0(x_1) \, \phi(x_1) \, x_1\,  dx_1
 + \int_{\Omega_0} \theta(x_2) \, u_0(x_1) \, \phi(x_1) \, x_1 \,dx_1 dx_2 \nonumber \\
= \int_0^1 \hat f(x_1) \, \phi(x_1) \, x_1 \,dx_1 
\end{eqnarray}
Due to $\xi^* \frac{\partial}{\partial x_1}(\phi \, x_1) = \xi^* x_1 \frac{\partial \phi}{\partial x_1} + \xi^*  \phi$, we can rewrite (\ref{limitP1}) as 
\begin{eqnarray} \label{limitP3}
\int_{\Omega_0} \Big\{ \xi^* \frac{\partial}{\partial x_1}(\phi x_1) - \phi \xi^* -    q\frac{\partial \phi}{\partial x_1} u_0
\Big\}  dx_1 dx_2  + \int_0^1 p \, u_0 \, \phi\, x_1 \, dx_1
+ \int_{\Omega_0} \theta \, u_0 \, \phi \, x_1\,dx_1 dx_2\nonumber \\
= \int_0^1 \hat f\, \phi \, x_1 \,dx_1 \; \forall \phi \in \mathcal{C}^\infty_0(0,1).
\end{eqnarray}
Therefore, it follows from (\ref{limitP2}) and (\ref{limitP3}) that, for all $\phi \in \mathcal{C}^\infty_0(0,1)$ 
\begin{equation}\label{xi}
0= \int_{\Omega_0} \Big\{  \phi \xi^* + q\frac{\partial \phi}{\partial x_1} u_0\Big\}  dx_1 dx_2= \int_{\Omega_0} \Big\{  \phi \xi^* - q\frac{\partial u_0}{\partial x_1} \phi \Big\}  dx_1 dx_2
\end{equation}
With the definition of $\hat q$ given by \eqref{definition-q-hat} 
and performing an iterated integration in (\ref{xi}) we obtain
$$
\int_0^1 \phi(x_1) \Big( \int_{-h_0}^{g_1} \xi^*(x_1, x_2) dx_2 - \hat q \frac{\partial u_0}{\partial x_1} \Big) dx_1=0 \qquad \forall \, \phi\, \in\, \mathcal{C}^\infty_0(0,1)
$$
So, the equation satisfied by $\xi^*$ is:
$$
 \int_{-h_0}^{g_1} \xi^*(x_1, x_2) dx_2= \hat q \frac{\partial u_0}{\partial x_1}
$$

The last step is placing this last equality in (\ref{limitP}). We get
\begin{equation} \label{EQfin}
\int_{0}^{1} \Big\{\hat q \frac{\partial u_0}{\partial x_1}\frac{\partial \varphi}{\partial x_1} 
+ \frac{|Y^*|}{L_1} \, u_0 \, \varphi + u_0\, \varphi \, p \,\Big\} dx_1 = \int_0^1  \hat f \, \varphi \, dx_1,\quad \forall \varphi\in H^1(0,1).
\end{equation}

Hence $u_0$ is the unique solution of (\ref{EQfin}), and we obtain that any convergent subsequence of $\{u^{\eps}\}$ tends to this unique solution. This complete the proof of Theorem \ref{teorema}.
\qedsymbol

\end{document}